\documentclass[a4paper,10.5pt]{article}
\usepackage[left=3cm,right=3cm,top=2.5cm,bottom=2.5cm]{geometry}

\usepackage{graphicx}
\usepackage{siunitx}
\usepackage{pgfplots}
\usepackage{bm}
\usepackage{nicefrac}
\usepackage{amssymb}
\usepackage{amsmath}
\usepackage[font+=small]{caption}
\usepackage[font+=small]{subcaption}
\usepackage{color,soul}%
\usepackage{tikz}
\usetikzlibrary{arrows,shapes,positioning}
\usetikzlibrary{plotmarks}
\usepackage{placeins}
\usepackage{booktabs}
\usepackage[noblocks,affil-it]{authblk}

\newcommand{\rmc}{{\mathrm{c}}}
\newcommand{\rme}{{\mathrm{e}}}
\newcommand{\rmf}{{\mathrm{f}}}
\newcommand{\rmt}{{\mathrm{t}}}
\newcommand{\rms}{{\mathrm{s}}}
\newcommand{\rmu}{{\mathrm{u}}}

\newcommand{\rmp}{{\mathrm{p}}}
\newcommand{\rmv}{{\mathrm{v}}}
\newcommand{\rmd}{{\mathrm{d}}}

\title{On Mixed Isogeometric Analysis of Poroelasticity}

\author[1]{Yared W. Bekele\thanks{Corresponding author. \textit{Email addresses}: yared.bekele@sintef.no (Yared W. Bekele), eivind.fonn@sintef.no (Eivind Fonn), trond.kvamsdal@math.ntnu.no (Trond Kvamsdal), arne.morten.kvarving@sintef.no (Arne M. Kvarving), steinar.nordal@ntnu.no (Steinar Nordal)}}
\author[2]{Eivind Fonn}
\author[2,3]{Trond Kvamsdal}
\author[2]{Arne M. Kvarving}
\author[4]{Steinar Nordal}
\affil[1]{{\scriptsize Rock and Soil Mechanics Group, SINTEF Building and Infrastructure, Trondheim, Norway}}
\affil[2]{{\scriptsize Department of Applied Mathematics, SINTEF ICT, Trondheim, Norway}}
\affil[3]{{\scriptsize Department of Mathematical Sciences, NTNU, Trondheim, Norway}}
\affil[4]{{\scriptsize Department of Civil and Environmental Engineering, NTNU, Trondheim, Norway}}

\date{}

\begin{document}

\maketitle

\begin{abstract}
    Pressure oscillations at small time steps have been known to be an issue in poroelasticity simulations. A review of proposed approaches to overcome this problem is presented. Critical time steps are specified to alleviate this in finite element analyses. We present a mixed isogeometric formulation here with a view to assessing the results at very small time steps. Numerical studies are performed on Terzaghi's problem and consolidation of a layered porous medium with a very low permeability layer for varying polynomial degrees, continuities across knot spans and spatial discretizations. Comparisons are made with equal order simulations.
\end{abstract}

\section{Introduction}
\label{sec:intro}

The study of porous materials, where the flow of fluid and solid deformation are
coupled, is essential in several areas of science and engineering. The theory of
poroelasticity is a mathematical formulation developed to describe these coupled
processes and predict the response of fluid saturated/unsaturated porous media
to external loading. There are different types of porous materials that are
studied under this theory such as soil, rock, concrete and other man-made
materials. Poroelasticity has a wide range of applications in different
disciplines of engineering mechanics and natural sciences. Some of the
application areas include geomechanics, biomechanics, reservoir engineering and
earthquake engineering. In addition to these diverse areas of application, it is
gaining popularity in the study of modern man-made porous media in material
science.

The mathematical formulations describing the fluid-solid coupled processes are
developed based on porous media theory where the multiphase medium is
approximated as a continuum, \cite{de1988historical}. The volume fraction
concept is used for averaging the properties of the multiphase medium in a
continuum formulation.

The governing partial differential equations of poroelasticity were first
developed for a one-dimensional case by
Terzaghi~\cite{terzaghi1923berechnung,terzaghi1925erdbaumechanik}. The
formulations were later generalized for a three-dimensional case and extended by
Biot~\cite{biot1941general,biot1955theory,biot1956general}. The mathematical
formulations have been studied extensively by several researchers since then.
Various analytical and numerical studies have been proposed in the literature.
Analytical solutions were obtained for problems with simplified material domains
and boundary conditions. Application to boundary value problems with complex
material domains and boundary conditions required the use of numerical methods.
The emergence of the finite element method opened the door for a detailed
numerical study of poroelasticity and for application to arbitrary geometries
and boundary conditions.

The finite element method was first applied to the governing equations of
poroelasticity to solve the initial boundary value problem of flow in a
saturated porous elastic medium by Sandhu and Wilson~\cite{sandhu1969finite}.
Hwang \textit{et al.}~\cite{hwang1971solutions} also used the finite element
method for plane strain consolidation problems and verified the results against
closed form solutions. The application of the finite element method started
gaining momentum afterwards and several researchers engaged themselves not only
on application problems but also in the investigation of the numerical
properties of the method within the context of poroelasticity. Ghaboussi and
Wilson~\cite{ghaboussi1973flow} applied the finite element method to partially
saturated elastic porous media and first noticed the ill-conditioning of the
matrix equations that may result when an incompressible fluid is assumed to
occupy the pore spaces. Booker and Small~\cite{booker1975investigation}
investigated the stability of the numerical solution when the finite element
method is applied to Biot's consolidation equations. The stability was studied
for different numerical integration schemes and time-step sizes. The numerical
performance of some finite element schemes for analysis of seepage in porous
elastic media was studied by Sandhu \textit{et al.}~\cite{sandhu1977numerical}.
They studied various spatial and temporal discretization schemes and evaluated
the numerical performances against the analytical solution of Terzaghi's
one-dimensional consolidation problem. Triangular and quadrilateral elements
with equal and mixed orders of interpolation for the displacement and pressure
were considered. It was shown that the elements with equal orders of
interpolation showed oscillatory behavior in the solution. Vermeer and
Verruijt~\cite{vermeer1981accuracy} derived a lower bound for the time-step size
in the analysis of consolidation by finite elements in terms of the mesh size
and the coefficient of consolidation. They showed that there is an accuracy
condition in the finite element analysis of consolidation by using a critical
time-step, below which oscillatory solutions are observed. The derived critical
time-step is strictly valid for a one-dimensional case and a uniform finite
element mesh. Reed~\cite{reed1984investigation} analyzed the numerical errors in
the analysis of consolidation by finite elements. It was shown that the use of a
mixed formulation for the field variables helps in reducing the pore pressure
oscillations but may not remove them entirely. They instead used Gauss point
smoothing to eliminate the pore pressure oscillations. Special finite elements
for the analysis of consolidation were proposed by Sandhu \textit{et
    al.}~\cite{sandhu1985special}. They presented ``singularity'' elements to model
pore pressures in the vicinity of free-draining loaded surfaces immediately
after application of loads. The elements were special in that they use special
interpolation schemes which reflect the actual variation of the field variables.

The finite element method became a well-established method for the analysis of
poroelasticity problems and the mathematical properties of the governing
equations and the numerical solution were studied in a further great detail.
Murad and Loula~\cite{murad1992improved} presented numerical analysis and error
estimates of finite element approximations of Biot's consolidation problem. They
used a mixed formulation and improved the rates of convergence by using a
sequential Galerkin Petrov-Galerkin post-processing technique. In a further
study,~\cite{murad1994stability}, they investigated the stability and
convergence of finite elements approximations of poroelasticity. They derived
decay functions showing that the pore pressure oscillations, arising from an
unstable approximation of the incompressibility constraint on the initial
conditions, decay in time. Finite element analysis of consolidation with
automatic time-stepping and error control was presented by Sloan and
Abbo~\cite{sloan1999biot1,sloan1999biot2}. Automatic time increments were
selected such that the temporal discretization error in the displacements is
close to a specified tolerance. Ferronato \textit{et
    al.}~\cite{ferronato2001ill} studied the ill-conditioning of finite element
poroelasticity equations with a focus on the instabilities that may affect the
pore pressure solution. They claim that the origin of most instabilities is due
to the assumption that, for initial conditions, the porous medium behaves as an
incompressible medium if the pore fluid is incompressible. They also argue that
oscillatory pore pressure solutions may not always be observed for very stiff
and low permeable materials depending on the critical time step. Gambolati
\textit{et al.}~\cite{gambolati2001numerical} studied the numerical performance
of projection methods in finite element consolidation models. Dureisseix
\textit{et al.}~\cite{dureisseix2003latin} proposed a large time increment
(LATIN) computational strategy for problems of poroelasticity to improve the
efficiency of the finite element analysis. A finite element formulation to
overcome spatial pore pressure oscillations caused by small time increments was
proposed by Zhu \textit{et al.}~\cite{zhu2004numerical}. Korsawe \textit{et
    al.}~\cite{korsawe2006finite} compared standard and mixed finite element methods
for poroelasticity. In particular, Galerkin and least-squares mixed finite
element methods were compared. They claim that Galerkin's method is able to
preserve steep pressure gradients but overestimates the effective stresses. On
the other hand, a least-squares mixed method was noticed to have the advantage
of direct approximation of the primary variables and explicit approximation of
Neumann type boundary conditions but to be computationally more expensive. A
mixed least-squares finite element method for poroelasticity was also proposed
by Tchonkova \textit{et al.}~\cite{tchonkova2008new}, claiming that pore
pressure oscillations are eliminated for different temporal discretizations. A
coupling of mixed and continuous Galerkin finite element methods for
poroelasticity was investigated for continuous and discrete in time cases by
Phillips and Wheeler~\cite{phillips2007coupling1,phillips2007coupling2}. They
also studied a coupling of mixed and discontinuous Galerkin finite-element
methods, \cite{phillips2008coupling3}. Haga \textit{et al.}~\cite{haga2012causes} studied the causes of pressure oscillations in low-permeable and low-compressible media by presenting two, three and four field mixed formulations in terms of the field variables displacement, pore fluid pressure, fluid velocity and solid skeleton stress.

A posteriori error estimation and adaptive refinement in poroelasticity has been studied by very few researchers. Larsson and Runesson~\cite{larsson2015sequential} presented a novel approach for space-time adaptive finite element analysis for the coupled consolidation problem in geomechanics. El-Hamalawi and Bolton~\cite{hamalawi1999posteriori} proposed an a posteriori error estimator for plane-strain geotechnical analyses based on superconvergent patch recovery with application to Biot's consolidation problem. They later extended the application of the a posteriori estimator for axisymmetric geotechnical analyses in~\cite{hamalawi2002posteriori}. Adaptive isogeometric finite element analysis, with LR B-Splines, for steady-state groundwater flow problems was presented by Bekele \textit{et al.}~\cite{bekele2015adaptive} but application to poroelasticity still remains as a task to study.

Isogeometric finite element analysis of poroelasticity was first presented by
Irzal \textit{et al.}~\cite{irzal2013isogeometric}. The advantages of the
smoothness of the basis functions in isogeometric analysis were highlighted in
their application. One of the advantages of higher continuity is that the
numerical implementation results in a locally mass conserving flow between
knotspans, analogous to elements in finite element analysis. But the formulation
presented relied on equal orders of interpolation for the field variables in
poroelasticity, namely displacement and pore fluid pressure. Such a formulation,
while still useful for several applications without significant numerical
challenges, has limitations when it comes to problems where the material
properties or boundary conditions are problematic.

In this paper, we present a mixed isogeometric formulation for poroelasticity.
To our best knowledge, this is the first time that a mixed isogeometric
formulation for poroelasticity is presented. The paper is structured as follows.
In Section~\ref{sec:equations}, the governing equations of poroelasticity are
presented. The fundamentals of isogeometric analysis and its particular features
of interest within the current context are discussed in Section~\ref{sec:iga}.
Numerical examples are given in Section~\ref{sec:numerical} and the observations
are summarized with concluding remarks in Section~\ref{sec:conclusions}.

\section{Governing Equations}
\label{sec:equations}

Biot's poroelasticity theory~\cite{biot1941general,biot1955theory} couples elastic solid deformation with fluid flow in the porous medium where the fluid flow is assumed to be governed by Darcy's law. The governing equations of the theory, the necessary boundary conditions, weak formulation and Galerkin finite element discretization are presented in the following sections.

\subsection{Linear Momentum Balance Equation}

The linear momentum balance equation for a fluid-saturated porous medium is given by:
\begin{equation}
\nabla \cdot \underbrace{
    \left( \bm\sigma' + \alpha p^\rmf \bm I \right)
}_{=\bm\sigma}
+ \rho \bm b
= \bm 0
\end{equation}  
where $\bm\sigma$ is the total stress, $\bm\sigma'$ is the effective stress, $\alpha$ is Biot's coefficient, $p^\rmf$ is the fluid pressure, $ \bm I $ is an identity matrix, $\rho$ is the overall density of the porous medium and $\bm b$ represents body forces. The Biot coefficient $\alpha$ can be calculated from:
\begin{equation}
\alpha = 1 - \frac{K_\rmt}{K_\rms}
\end{equation}
where $K_\rmt$ and $K_\rms$ are the bulk moduli of the porous medium and solid
particles, respectively.

The constitutive equation for poroelasticity relates stress and strain linearly as:
\begin{equation}
\bm\sigma' = \bm D : \bm\varepsilon
\end{equation}
where $\bm D$ is a fourth-order stiffness tensor. Small deformations are also
assumed, so the strain $\bm\varepsilon$ satisfies a linear first-order equation with respect to the displacement $\bm u$,
\begin{equation}
\bm\varepsilon = \frac{1}{2} \left( \nabla\bm u + \nabla^\intercal \bm u \right)
\end{equation}
where $\nicefrac{1}{2}(\nabla+\nabla^\intercal)$ is the symmetrized gradient
operator i.e.
\begin{equation}
\varepsilon_{ij} = \frac{1}{2} \left(
\frac{\partial u_i}{\partial x_j} + \frac{\partial u_j}{\partial x_i}
\right).
\end{equation}
In the following, it will be convenient to lower tensors and higher differential operators to Voigt notation, which represents the symmetric $d \times d$ tensor $\bm \sigma'$ as a $\nicefrac{d(d+1)}{2}$-vector, which we will denote with a tilde:
\[
\underbrace{
    \left[ \begin{matrix}
    \sigma'_{xx} & \sigma'_{xy} & \sigma'_{xz} \\
    \sigma'_{xy} & \sigma'_{yy} & \sigma'_{yz} \\
    \sigma'_{xz} & \sigma'_{yz} & \sigma'_{zz} \\
    \end{matrix} \right] 
}_{\bm\sigma'}
\Longleftrightarrow
\underbrace{
    \left\lbrace \begin{matrix}
    \sigma'_{xx} & \sigma'_{yy} & \sigma'_{zz} &
    \sigma'_{yz} & \sigma'_{xz} & \sigma'_{xy}
    \end{matrix} \right\rbrace ^\intercal
}_{\tilde{\bm\sigma}'}
\]
A similar conversion takes place for the strains, where the shear strains are replaced by the \emph{engineering} shear strains:
\[
\underbrace{
    \left[ \begin{matrix}
    \varepsilon_{xx} & \varepsilon_{xy} & \varepsilon_{xz} \\
    \varepsilon_{xy} & \varepsilon_{yy} & \varepsilon_{yz} \\
    \varepsilon_{xz} & \varepsilon_{yz} & \varepsilon_{zz} \\
    \end{matrix}\right] 
}_{\bm\varepsilon}
\Longleftrightarrow
\underbrace{
    \left\lbrace \begin{matrix}
    \varepsilon_{xx} & \varepsilon_{yy} & \varepsilon_{zz} &
    2\varepsilon_{yz} & 2\varepsilon_{xz} & 2\varepsilon_{xy}
    \end{matrix}\right\rbrace ^\intercal
}_{\tilde{\bm\varepsilon}}
\]
Voigt notation allows us to express the equilibrium equation and the stress-strain equation using the same differential operator $\bm L$,
\begin{equation}
{\bm L}^\intercal = \left[ 
\begin{matrix}
\frac{\partial}{\partial x} & 0 & 0 & \frac{\partial}{\partial y} & 0 & \frac{\partial}{\partial z} \\
0 & \frac{\partial}{\partial y} & 0 & \frac{\partial}{\partial x} & \frac{\partial}{\partial z} & 0 \\
0 & 0 & \frac{\partial}{\partial z} & 0 & \frac{\partial}{\partial y} & \frac{\partial}{\partial x}
\end{matrix}
\right] 
\end{equation}
Using $ \bm L $ yields the following equilibrium equation in terms of the two primary unknowns $\bm u$ and $p^\rmf$,
\begin{equation} \label{eqn:equilibrium}
\bm L^\intercal \tilde{\bm D} \bm L \bm u
- \alpha \nabla p^\rmf
+ \rho \bm b
= \bm 0,
\end{equation}
where $\tilde{\bm D}$ is the Voigt notation equivalent of $\bm D$, taking into account the aforementioned engineering shear strains. We will generally assume isotropic materials, where $\tilde{\bm D}$ takes the block form (in terms of Young's modulus $E$ and Poisson's ratio $\nu$)
\begin{equation}
\tilde{\bm D} =
\frac{E}{(1 + \nu)(1 - 2\nu)}
\left[ \begin{matrix}
\tilde{\bm D}_{11} & \bm 0 \\
\bm 0 & \tilde{\bm D}_{22}
\end{matrix} \right] 
\end{equation}
where the two blocks are given as
\begin{align}
\begin{split}
\tilde{\bm D}_{11} &= (1-2\nu) \bm I + \nu \boldsymbol{\mathit{1}} \\
\tilde{\bm D}_{22} &= \frac{1-2\nu}{2} \bm I
\end{split}
\end{align}  
and $ \boldsymbol{\mathit{1}} $ is a matrix of ones.

\subsection{Mass Balance Equation}

A mass conservation equation together with the equilibrium equation in~\eqref{eqn:equilibrium} completes the governing equations of poroelasticity. The \emph{fluid content} $\zeta$ is given by
\begin{equation}
\zeta = \alpha \nabla \cdot \bm u + c p^\rmf
\end{equation}
where $c$ is the \emph{storativity} or \emph{specific storage coefficient} at
constant strain. It is given by
\begin{equation}
c = \frac{\alpha - n}{K_\rms} + \frac{n}{K_\rmf}
\end{equation}
where $K_\rmf$ is the bulk modulus of the fluid and $n$ is the porosity of the material. The change in the fluid content $\zeta$ satisfies the equation
\begin{equation}
\frac{\partial \zeta}{\partial t} + \nabla \cdot \bm w = 0
\end{equation}
where $\bm w$ is the fluid flux, which is given by Darcy's law as:
\begin{equation}
\bm w = -\frac{1}{\gamma_\rmf} \bm k \cdot \left(\nabla p^\rmf - \rho_\rmf \bm b \right)
\end{equation}
where $ \gamma_\rmf $ is the unit weight of the fluid, $ \rho_\rmf $ its density and $ \bm k $ is the hydraulic conductivity matrix.

The final equation of mass balance is then
\begin{equation} \label{eqn:mass-balance}
\alpha \nabla \cdot \dot{\bm u}
+ c \frac{\partial p^\rmf}{\partial t} 
+ \nabla \cdot \left[ -\frac{1}{\gamma_\rmf} \bm k \cdot \left(\nabla p^\rmf - \rho_\rmf \bm b \right) \right]  
= 0.
\end{equation}

\subsection{Boundary Conditions}

The governing linear momentum and mass balance equations in~\eqref{eqn:equilibrium} and \eqref{eqn:mass-balance}, respectively, are accompanied by the usual boundary conditions in the formulation of bounary value problems. Let $(\Gamma_D^\rmu,\Gamma_D^\rmp)$ and
$(\Gamma_N^\rmu,\Gamma_N^\rmp)$ be two partitions of the boundary $\partial \Omega$ of domain $ \Omega $, for representing Dirichlet and Neumann boundary conditions, respectively.

The Dirichlet boundary conditions for the equilibrium~\eqref{eqn:equilibrium} and mass balance~\eqref{eqn:mass-balance} equations are
\begin{equation}
\begin{cases}	     
\bm u = \overline{\bm u} \qquad & \text{on $\Gamma_D^\rmu$}, \\
p^\rmf = \overline{p}^\rmf \qquad & \text{on $\Gamma_D^\rmp$},
\end{cases}
\end{equation}
where $\overline{\bm u}$ and $\overline{p}_\rmf$ are the prescribed displacement and pressure, respectively.

The Neumann boundary conditions are
\begin{equation}
\begin{cases}
\bm \sigma \cdot \bm n = \overline{\bm t} \qquad & \text{on $\Gamma_N^\rmu$}, \\
\bm w \cdot \bm n = \overline{q} \qquad & \text{on $\Gamma_N^\rmp$},
\end{cases}
\end{equation}
where $\bm n$ is the outward pointing normal vector, $\overline{\bm t}$ is the surface traction and $\overline{q}$ is the fluid flux on the boundary.

\subsection{Variational Formulation}

To derive the variational formulations of equations
\eqref{eqn:equilibrium} and \eqref{eqn:mass-balance}, we introduce a vector-valued test function $\delta \bm u$, which vanishes on $\Gamma_D^\rmu$, and a scalar test function $\delta p$, which vanishes on $\Gamma_D^\rmp$.

We start with the total stress formulation of the linear momentum balance equation, which from equation~\eqref{eqn:equilibrium} is given by
\begin{equation}
\nabla \cdot \bm \sigma + \rho \bm b = 0.
\end{equation}
Multiplying by the test function $\delta \bm u$ and integrating over the domain $ \Omega $ gives
\begin{equation}
\int_\Omega \delta \bm u^\intercal \nabla \cdot \bm \sigma \rmd \Omega + \int_\Omega \delta \bm u^\intercal \rho \bm b \rmd \Omega = \bm 0. 
\end{equation}
The first term in the above equation contains a double derivative of the unknown displacement, and is relaxed using a form of Green's theorem,
\begin{align}
\begin{split}
\int_\Omega \delta \bm u^\intercal \nabla \cdot \bm \sigma \rmd \Omega
&= \sum_i \int_\Omega \delta u_i \nabla \cdot \bm \sigma_i \rmd \Omega \\
&= \sum_i \int_{\partial \Omega} \delta u_i \bm \sigma_i \cdot \bm n \rmd \Gamma
- \sum_i \int_\Omega \bm \nabla \delta u_i \cdot \sigma_i \rmd \Omega \\
&= \int_{\Gamma_N^\rmu} \delta \bm u^\intercal \overline{\bm t} \rmd \Gamma
- \int_\Omega \nabla \delta \bm u : \bm \sigma \rmd \Omega.
\end{split}
\end{align}
Due to the symmetry of the stress tensor, the last term is expressible in Voigt notation,
\begin{equation}
\nabla \delta \bm u : \bm \sigma = {(\bm L \delta \bm u)}^\intercal \bm \sigma,
\end{equation}
yielding the weak form of~\eqref{eqn:equilibrium} as
\begin{equation} \label{eqn:weak:equilibrium}
\int_\Omega {(\bm L \delta \bm u)}^\intercal \tilde{\bm D} (\bm L \bm u) \rmd \Omega
- \alpha \int_\Omega {(\bm L \delta \bm u)}^\intercal \tilde{\bm I} p^\rmf \rmd \Omega
= \int_\Omega \delta \bm u^\intercal \rho \bm b \rmd \Omega
+ \int_{\Gamma_N^\rmu} \delta \bm u^\intercal \overline{\bm t} \rmd \Gamma
\end{equation}
where we have used $\tilde{\bm I}$ as the Voigt notation identity operator, which for a general three-dimensional case is given by:
\begin{equation}
\tilde{\bm I} = \left\lbrace 1,1,1,0,0,0 \right\rbrace^\intercal
\end{equation}
For the mass balance equation, multiplying~\eqref{eqn:mass-balance} by the scalar test function $ \delta p $ and integrating over the domain $ \Omega $, we get
\begin{equation}
\alpha \int_\Omega \delta p \nabla \cdot \dot{\bm u} \rmd \Omega
+ c \int_\Omega \delta p \frac{\partial p^\rmf}{\partial t} \rmd \Omega
+ \int_\Omega \delta p \nabla \cdot \left[ -\frac{1}{\gamma_\rmf} \bm k \cdot \left(\nabla p^\rmf - \rho_\rmf \bm b \right) \right] \rmd \Omega
= 0.
\end{equation}
Again, by applying Green's theorem to the last term, we obtain
\begin{equation}
\int_\Omega \delta p \nabla \cdot \left[ -\frac{1}{\gamma_\rmf} \bm k \cdot \left(\nabla p^\rmf - \rho_\rmf \bm b \right) \right] \rmd \Omega
= \int_{\Gamma_N^\rmp} \delta p \overline{q} \rmd \Gamma
- \int_\Omega \nabla \delta p \cdot \left[ -\frac{1}{\gamma_\rmf} \bm k \cdot \left(\nabla p^\rmf - \rho_\rmf \bm b \right) \right] \rmd \Omega.
\end{equation}
Thus, the weak form of the mass balance equation,~\eqref{eqn:mass-balance}, is
\begin{equation} \label{eqn:weak:mass-balance}
\alpha \int_\Omega \delta p \nabla \cdot \dot{\bm u} \rmd \Omega
+ c \int_\Omega \delta p \frac{\partial p^\rmf}{\partial t} \rmd \Omega
+ \int_\Omega \nabla \delta p^\intercal \frac{1}{\gamma_\rmf} \bm k \nabla p^\rmf \rmd \Omega
= \int_\Omega \nabla \delta p^\intercal \frac{1}{\gamma_\rmf} \bm k \rho_\rmf \bm b \rmd \Omega
- \int_{\Gamma_N^\rmp} \delta p \overline{q} \rmd \Gamma.
\end{equation}
\subsection{Galerkin Finite Element Formulation}

With a suitable number $N$ of basis functions defined, let $\bm N_\rmp : \Omega \to \mathbb{R}^{1 \times N}$ and $\bm N_\rmu : \Omega \to \mathbb{R}^{d \times dN}$
be the basis interpolation matrices for the pressure and displacement
respectively. The unknowns and the test functions can then be represented using coefficient vectors:
\begin{align}  
\begin{aligned}
\bm u &= \bm N_\rmu \bm u^\rmc, &\delta \bm u &= \bm N_\rmu \delta \bm u^\rmc, \\
p^\rmf &= \bm N_\rmp \bm p^\rmc, & \delta p &= \bm N_\rmp \delta \bm p^\rmc 
\end{aligned}
\label{eqn:galerkin_coeff}
\end{align}
where $ \bm u^\rmc $ and $ \bm p^\rmc $ are the control point values of the displacement and pressure field variables. Application of~\eqref{eqn:galerkin_coeff} to the weak form of the linear momentum balance equation in~\eqref{eqn:weak:equilibrium} results in the matrix the discrete system of equations~(after canceling $\delta \bm u^\rmc$ and $\bm \delta p^\rmc$, as equations~\eqref{eqn:equilibrium} and~\eqref{eqn:mass-balance} are supposed to be valid for any choice of these)
\begin{equation}
\bm K \bm u^\rmc - \bm Q \bm p^\rmc = \bm f_\rmu
\label{eqn:equilibrium_galerkin}
\end{equation}
where the stiffness matrix $ \bm K $, the coupling matrix $ \bm Q $ and the vector of body forces and surface tractions $ \bm f_u $ are given by
\begin{align}
\begin{split}
\bm K &= \int_\Omega \bm B^\intercal \tilde{\bm D} \bm B \rmd \Omega
\qquad \\
\bm Q &= \int_\Omega \bm B^\intercal \alpha \tilde{\bm I} \bm N_\rmp \rmd \Omega \\
\bm f_\rmu &= \int_\Omega \bm N_\rmu^\intercal \rho \bm b \rmd \Omega
+ \int_{\Gamma_N^\rmu} \bm N_\rmu^\intercal \overline{\bm t} \rmd \Gamma
\end{split}
\end{align}
Here $ \bm B = \bm L \bm N_\rmu $ is the strain-displacement matrix. Similarly, using~\eqref{eqn:galerkin_coeff} in the weak form of the mass balance equation in~\eqref{eqn:weak:mass-balance} results in the discrete system of equations
\begin{equation}
\bm Q^\intercal \frac{\partial \bm u^\rmc}{\partial t} + \bm S \frac{\partial \bm p^\rmc}{\partial t}  + \bm P \bm p^\rmc = \bm f_\rmp
\label{eqn:mass_galerkin}
\end{equation}
where the storage matrix $ \bm S $, the permeability matrix $ \bm P $ and the vector of fluid body forces and fluxes $ \bm f_\rmp $ are given by
\begin{align}
\begin{split}
\bm S &= \int_\Omega \bm N_\rmp^\intercal c \bm N_\rmp \rmd \Omega \\
\bm P &= \int_\Omega \nabla \bm N_\rmp^\intercal \frac{1}{\gamma_\rmf} \bm k \nabla \bm N_\rmp \rmd \Omega \\
\bm f_\rmp &= \int_\Omega \nabla \bm N_\rmp^\intercal \frac{1}{\gamma_\rmf} \bm k \rho_\rmf \bm b \rmd \Omega - \int_{\Gamma_N^\rmp} \bm N_\rmp^\intercal \overline{q} \rmd \Gamma.
\end{split}
\end{align}
Combining equations~\eqref{eqn:equilibrium_galerkin} and \eqref{eqn:mass_galerkin} results in the coupled system of equations for poroelasticity
\begin{equation} \label{eqn:discr-asymmetric}
\left[ \begin{matrix} \bm 0 & \bm 0 \\ \bm Q^\intercal & \bm S \end{matrix} \right] 
\left\lbrace \begin{matrix} \dot{\bm u}^\rmc \\ \dot{\bm p}^\rmc \end{matrix} \right\rbrace 
+
\left[ \begin{matrix} \bm K & - \bm Q \\ \bm 0 & \bm P \end{matrix} \right] 
\left\lbrace \begin{matrix} \bm u^\rmc \\ \bm p^\rmc \end{matrix} \right\rbrace 
=
\left\lbrace \begin{matrix} \bm f_\rmu \\ \bm f_\rmp \end{matrix} \right\rbrace. 
\end{equation}
A symmetric system of equations can be obtained by time-differentiating the first equation and multiplying one of the equations by $ -1 $, \cite{lewis1998finite}:
\begin{equation} \label{eqn:discr-symmetric}
\left[ \begin{matrix} - \bm K & \bm Q \\ \bm Q^\intercal & \bm S \end{matrix} \right] 
\left\lbrace \begin{matrix} \dot{\bm u}^\rmc \\ \dot{\bm p}^\rmc \end{matrix} \right\rbrace 
+
\left[ \begin{matrix} \bm 0 & \bm 0 \\ \bm 0 & \bm P \end{matrix} \right] 
\left\lbrace \begin{matrix} \bm u^\rmc \\ \bm p^\rmc \end{matrix} \right\rbrace 
=
\left\lbrace \begin{matrix} - \dot{\bm f}_\rmu \\ \bm f_\rmp \end{matrix} \right\rbrace   
\end{equation}
In this formulation, it is important that time-dependent quantities involved in $\bm f_u$, such as traction and body forces, are ``ramped up'' from an initial equilibrium instead of being applied immediately. This can be done in the first time step.

\subsection{Temporal Discretization}

The generalized trapezoidal rule (GTR) is applied for the temporal discretization of the coupled system of matrix equations in~\eqref{eqn:discr-symmetric}. Representing the vector of unknowns by $ \bm X = \left\lbrace \bm u^c,\bm p^c \right\rbrace^\intercal $, we have the GTR approximation
\begin{align}
\begin{split}
\left. \frac{\partial \bm X}{\partial t} \right\rvert_{n+\theta} &= \frac{{{\bm X}_{n+1}}-{{\bm X}_{n}}}{\Delta t} \\
{{\bm X}_{n+\theta }} &= (1-\theta ){{\bm X}_{n}}+\theta {{\bm X}_{n+1}}
\end{split}
\label{eqn:GTR}
\end{align}
where $ \theta $ is a time integration parameter which has limits $ 0 \leq \theta \leq 1 $ and $ n $ is a time step identifier. Adopting backward Euler time stepping ($ \theta=1 $) with time step $\Delta t$ and applying~\eqref{eqn:GTR} to~\eqref{eqn:discr-symmetric} we obtain the system of equations
\begin{equation} \label{eqn:discr-time}
\left[ \begin{matrix}
- \bm K & \bm Q \\ \bm Q^\intercal & \bm S + \Delta t \bm P
\end{matrix} \right] 
\left\lbrace \begin{matrix} \bm u^\rmc \\ \bm p^\rmc \end{matrix}\right\rbrace_{n+1}
=
\left[ \begin{matrix}
- \bm K & \bm Q \\ \bm Q^\intercal & \bm S
\end{matrix} \right] 
\left\lbrace \begin{matrix} \bm u^\rmc \\ \bm p^\rmc \end{matrix}\right\rbrace_n
+ \Delta t
\left\lbrace \begin{matrix} - \dot{\bm f}_u \\ \bm f_p \end{matrix}\right\rbrace_{n+1}
\end{equation}
which is a linear system in this case, for poroelasticity, as the coefficient matrices are independent of the unknowns.

\section{Isogeometric Analysis}
\label{sec:iga}

\subsection{Introduction}

Since its first introduction by Hughes \textit{et al}. \cite{hughes2005iga}, isogeometric analysis (IGA) has been successfully applied to several areas of engineering mechanics problems. The fundamental aim for the introduction of IGA was the idea of bridging the gap between finite element analysis (FEA) and computer-aided design (CAD). The main concept behind the method is the application of the same basis functions used in CAD for performing finite element analysis. In the process of its application to various engineering problems, IGA has shown advantages over the conventional finite element method, for instance the ease of performing finite element analysis using higher order
polynomials.

We briefly present the fundamentals behind B-Splines and Non-Uniform Rational B-Splines (NURBS) in the next section and highlight the features of IGA that are important in our context.

\subsection{Fundamentals on B-Splines and NURBS}

We start the discussion on B-Splines and NURBS by first defining a \textit{knot
    vector}. A knot vector in one dimension is a non-decreasing set of coordinates
in the parameter space, written $ {\Xi} = \{\xi_1,\xi_2,...,\xi_{n+p+1}\} $,
where $ \xi_i \in \mathbb{R} $ is the $ i^{th} $ knot, $ i $ is the knot index,
$ i = 1,2,...,n+p+1 $, $ p $ is the polynomial order, and $ n $ is the number of
basis functions. Knot vectors may be uniform or non-uniform depending on whether
the knots are equally spaced in the parameter space or not.

A univariate B-Spline curve is parametrized by a linear combination of $ n $
B-Spline basis functions, $ \{N_{i,p}\}_{i=1}^n $. The coefficients
corresponding to these functions, $ \{{\bm X}_{i}\}_{i=1}^n $, are referred to
as control points. The B-Spline basis functions are recursively defined starting
with piecewise constants ($p=0$):

\begin{equation}
N_{i,0}(\xi) = \begin{cases}
1 & \text{if } \xi_i \leq \xi < \xi_{i+1} \\
0 & \text{otherwise}
\end{cases}
\end{equation}

For higher-order polynomial degrees ($ p \geq 1 $), the basis functions are
defined by the Cox-de Boor recursion formula:

\begin{equation}
N_{i,p}(\xi) =
\frac{\xi - \xi_i}{\xi_{i+p}-\xi_i}N_{i,p-1}(\xi) +
\frac{\xi_{i+p+1} - \xi}{\xi_{i+p+1}-\xi_{i+1}}N_{i+1,p-1}(\xi)
\end{equation}

B-Spline geometries, curves, surfaces and solids, are constructed from a linear
combination of B-Spline basis functions. Given $ n $ basis functions $ N_{i,p} $
and corresponding control points $ {\bm P}_i \in \mathbb{R}^d, i=1,2,...,n $, a
piecewise polynomial B-Spline curve is given by:

\begin{equation}
{\bm C}(\xi) = \sum\limits_{i=1}^n N_{i,p}(\xi) {\bm P}_i
\end{equation}

Similarly, for a given control net $ {\bm P}_{i,j}, i=1,2,...,n,j=1,2,...,m, $
polynomial orders $ p $ and $ q $, and knot vectors $ {\Xi} =
\{\xi_1,\xi_2,...,\xi_{n+p+1}\} $, and $ \mathcal{H} =
\{\eta_1,\eta_2,...,\eta_{m+q+1}\} $, a tensor product B-Spline surface is
defined by:

\begin{equation}
{\bm S}(\xi,\eta) =
\sum\limits_{i=1}^n \sum\limits_{j=1}^m
N_{i,p}(\xi) M_{j,q}(\eta) {\bm P}_{i,j}
\end{equation}

B-Spline solids are defined in a similar way as B-Spline surfaces from tensor
products over a control lattice.

NURBS are built from B-Splines to represent a wide array of objects that cannot
be exactly represented by polynomials. A NURBS entity in $ \mathbb{R}^d $ is
obtained by projective transformation of a B-Spline entity in $ \mathbb{R}^{d+1}
$. The control points for the NURBS geometry are found by performing exactly the
same projective transformation to the control points of the B-Spline curve. A
detailed treatment of B-Splines and NURBS can be referred from Cottrell
\textit{et al.} \cite{cottrell2009iga}.


\subsection{Important Features in Current Context}

IGA has a number of advantages over FEA such as the ability to represent exact
geometries of structures or domains, non-negative basis functions and
isoparametric mapping at patch level. In the context of the current work, we
focus on the features of IGA that are especially important. These features are
improved continuity because of the smoothness of the basis functions and the
ability to perform simulations with high continuity and high regularity meshes.
We look closely into each here.

\subsubsection{Continuity}

One of the most distinctive and powerful features of IGA is 
that the basis functions will be $ C^{p-m} $ continuous across knotspans 
(analogous to elements in FEA), where $ p $ is the polynomial degree and $ m $
is the multiplicity of the knot. This means that the continuity across knotspans
can be controlled by the proper choice of $ p $ and $ m $. The continuity can be
decreased by repeating a knot - important to model non-smooth geometry features
or to facilitate the application of boundary conditions. For instance, quadratic
($ p=2 $) splines are $ C^1 $ continuous over non-repeated knots while quadratic
Lagrange finite element bases are only $ C^0 $ continuous. If we consider the
quartic ($ p=4 $) basis functions constructed from the open, non-uniform knot 
vector $ \Xi = \left\lbrace 0,0,0,0,0,1,2,2,3,3,3,4,4,4,4,5,5,5,5,5 
\right\rbrace $, we get different continuities across knotspans as shown in 
Figure~\ref{fig:Cp-m}.

\subsubsection{$k$-refinement}

IGA and FEA both allow $ h $- and $ p $-refinements i.e. increasing the number of knotspans by knot insertion (increasing the number of elements in FEA) and raising the polynomial order. The non-commutativity of knot insertion and polynomial order elevation results in a type of refinement
that is unique to IGA, called $ k $-refinement. This is achieved by performing polynomial order elevation followed by knot insertion. This results in a high continuity mesh with the least number of degrees of freedom i.e. high regularity.

\definecolor{col1}{RGB}{138,82,73}
\definecolor{col2}{RGB}{214,183,60}
\definecolor{col3}{RGB}{114,194,214}
\definecolor{col4}{RGB}{138,39,118}
\definecolor{col5}{RGB}{214,81,60}
\definecolor{col6}{RGB}{183,214,60}
\definecolor{col7}{RGB}{107,145,201}
\definecolor{col8}{RGB}{214,114,154}
\definecolor{col9}{RGB}{214,154,114}
\definecolor{col10}{RGB}{99,138,73}
\definecolor{col11}{RGB}{39,65,138}
\definecolor{col12}{RGB}{138,78,39}
\definecolor{col13}{RGB}{60,214,101}
\definecolor{col14}{RGB}{140,114,214}
\definecolor{col15}{RGB}{138,105,39}
\definecolor{col16}{RGB}{60,214,183}
\definecolor{col17}{RGB}{204,60,214}

\begin{figure}[t]
    \centering
    \begin{tikzpicture}
    \begin{axis}[
    xmin=0, xmax=5,
    ymin=0, ymax=1,
    width=0.8\textwidth,
    height=0.3\textwidth,
    xtick={0,1,2,3,4,5},
    xticklabels={0 0 0 0 0,1,2 2,3 3 3,4 4 4 4,5 5 5 5 5},
    ytick={1},
    yticklabels={1},
    axis lines=left,
    clip=false,
    ]
    \addplot[mark=none, very thick, col1]
    table[x index={0}, y index={1}]{data/basis/basis.csv};
    \addplot[mark=none, very thick, col2]
    table[x index={0}, y index={2}]{data/basis/basis.csv};
    \addplot[mark=none, very thick, col3]
    table[x index={0}, y index={3}]{data/basis/basis.csv};
    \addplot[mark=none, very thick, col4]
    table[x index={0}, y index={4}]{data/basis/basis.csv};
    \addplot[mark=none, very thick, col5]
    table[x index={0}, y index={5}]{data/basis/basis.csv};
    \addplot[mark=none, very thick, col6]
    table[x index={0}, y index={6}]{data/basis/basis.csv};
    \addplot[mark=none, very thick, col7]
    table[x index={0}, y index={7}]{data/basis/basis.csv};
    \addplot[mark=none, very thick, col8]
    table[x index={0}, y index={8}]{data/basis/basis.csv};
    \addplot[mark=none, very thick, col9]
    table[x index={0}, y index={9}]{data/basis/basis.csv};
    \addplot[mark=none, very thick, col10]
    table[x index={0}, y index={10}]{data/basis/basis.csv};
    \addplot[mark=none, very thick, col11]
    table[x index={0}, y index={11}]{data/basis/basis.csv};
    \addplot[mark=none, very thick, col12]
    table[x index={0}, y index={12}]{data/basis/basis.csv};
    \addplot[mark=none, very thick, col13]
    table[x index={0}, y index={13}]{data/basis/basis.csv};
    \addplot[mark=none, very thick, col14]
    table[x index={0}, y index={14}]{data/basis/basis.csv};
    \addplot[mark=none, very thick, col15]
    table[x index={0}, y index={15}]{data/basis/basis.csv};
    \node[anchor=south] (zero) at (axis cs: 0,1.2){$C^{-1}$};
    \node[anchor=south] (one) at (axis cs: 1,1.2){$C^{3}$};
    \node[anchor=south] (two) at (axis cs: 2,1.2){$C^{2}$};
    \node[anchor=south] (three) at (axis cs: 3,1.2){$C^{1}$};
    \node[anchor=south] (four) at (axis cs: 4,1.2){$C^{0}$};
    \node[anchor=south] (five) at (axis cs: 5,1.2){$C^{-1}$};
    \draw[->,-stealth] (zero) -- (axis cs: 0,1.05);
    \draw[->,-stealth] (one) -- (axis cs: 1,0.50);
    \draw[->,-stealth] (two) -- (axis cs: 2,0.65);
    \draw[->,-stealth] (three) -- (axis cs: 3,0.70);
    \draw[->,-stealth] (four) -- (axis cs: 4,1.08);
    \draw[->,-stealth] (five) -- (axis cs: 5,1.05);
    \end{axis}
    \end{tikzpicture}
    \caption{Different continuities across knotspans, after~\cite{cottrell2009iga}.}
    \label{fig:Cp-m}
\end{figure}
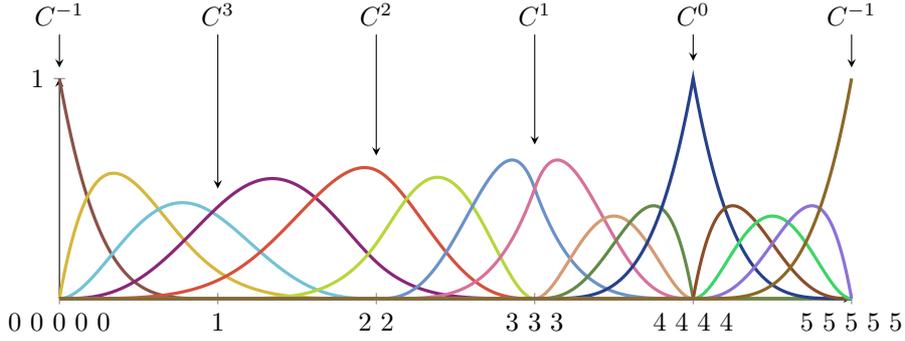

\FloatBarrier

\subsection{Mixed Isogeometric Formulation}

\begin{figure}
    \begin{center}
        \begin{tikzpicture}[scale=0.8]
        \draw [step=1,gray,thin] (0,0) grid (5,5);
        \fill[green!60!white] (2,2) rectangle (3,3);
        \foreach \a in {2,3} {
            \foreach \b in {2,3} {
                \node [circle,fill=blue,scale=0.5] at (\b,\a) {};
            }
        }
        \node [below] at (2.5,0) {$ p=1 $};
        \draw [step=1,gray,thin] (6,0) grid (11,5);
        \fill[green!60!white] (8,2) rectangle (9,3);
        \foreach \a in {1.5,2.5,3.5} {
            \foreach \b in {7.5,8.5,9.5} {
                \node [circle,fill=blue,scale=0.5] at (\b,\a) {};
            }
        }
        \node [below] at (8.5,0) {$ p=2 $};
        \draw [step=1,gray,thin] (12,0) grid (17,5);
        \fill[green!60!white] (14,2) rectangle (15,3);
        \foreach \a in {1,2,3,4} {
            \foreach \b in {13,14,15,16} {
                \node [circle,fill=blue,scale=0.5] at (\b,\a) {};
            }
        }
        \node [below] at (14.5,0) {$ p=3 $};  
        \end{tikzpicture}
    \end{center}
    \caption{Number of control points for a given element on a simple B-Spline surface with different polynomial degrees. The element is highlighted and the blue squares represent control points.}
    \label{fig:controlpoints}
\end{figure}
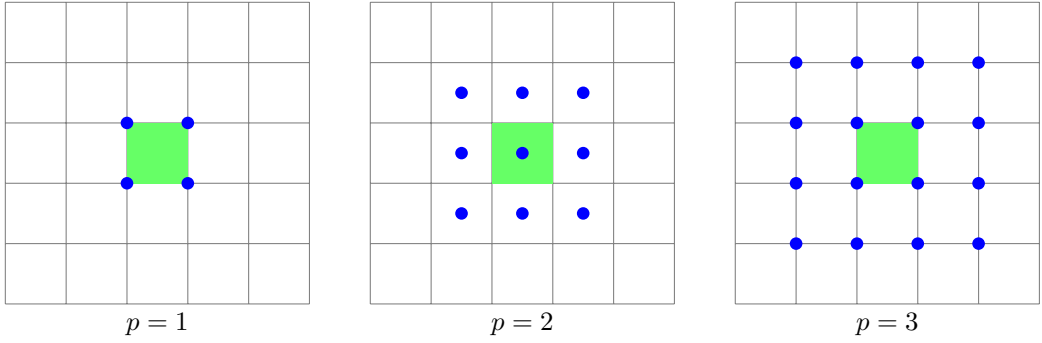

A mixed formulation is constructed by first defining the knot vectors and basis functions defining the geometry of the domain. The polynomial order defining the geometry is used as the polynomial degree for one of the field variables and is raised by the desired degree for the other field variable. In our context, the polynomial order for the pressure, $ p_\rmp $, is defined by the geometry construction and the polynomial order for the displacement, $ p_\rmu $, is raised by one. Both $ p_\rmp $ and $ p_\rmu $ can then be raised to the desired degree starting from the initial definition. For example, a simple two-dimensional geometry defined by the knot vectors $ \Xi = \left\lbrace 0,0,1,1 \right\rbrace $ and $ \mathcal{H} = \left\lbrace 0,0,1,1 \right\rbrace $ implies $ p_\rmp = 1 $ and $ p_\rmu = 2 $ with 4 and 9 control points, respectively.  The number of control points, location of degrees of freedom in IGA, on a B-Spline surface for different polynomial degrees is shown in Figure~\ref{fig:controlpoints}.

\FloatBarrier

\section{Numerical Examples}
\label{sec:numerical}

In this section, the performance of a mixed isogeometric formulation is investigated for some numerical examples. We first consider Terzaghi's classical one-dimensional consolidation problem for verification and mesh convergence studies. Consolidation of a layered medium with a low permeability layer sandwiched between two high permeability layers is studied. The mixed formulation results are compared with equal order simulation.

\subsection{Terzaghi's Problem}

\begin{figure}[h]
    \centering
    \begin{tikzpicture}[scale=4]
    \filldraw[fill=gray!15, draw=black] (0,0) rectangle (1,1);
    \foreach \x in {0, 0.2, 0.4, 0.6, 0.8, 1}
    \draw[thick, ->,-stealth] (\x,1.15) -- (\x,1);
    \draw [thick] (0,1.15) -- (1,1.15);
    \draw[<->,-stealth, dotted] (0.3,0) -- (0.3,1);
    \node[anchor=west] at (0.3, 0.5) {$h = \SI{8}{\milli\meter}$};
    \node[anchor=south] at (0.5, 1.15) {
        $\overline{t}_y = -p_0$,
        $\overline{u}_x = 0$,
        $\overline{p}^\rmf = 0$
    };
    \node[anchor=north] at (0.5, 0.0) {
        $\overline{\bm u} = \bm 0$,
        $\overline{q} = 0$
    };
    \node[anchor=west] at (1.0, 0.5) {
        $\begin{aligned}
        \overline{u}_x &= 0 \\
        \overline{q} &= 0
        \end{aligned}$
    };
    \node[anchor=east] at (0.0, 0.5) {
        $\begin{aligned}
        \overline{u}_x &= 0 \\
        \overline{q} &= 0
        \end{aligned}$
    };
    \end{tikzpicture}
    \caption{Terzaghi's problem: Domain and boundary conditions.}
    \label{fig:terzhagi}
\end{figure}
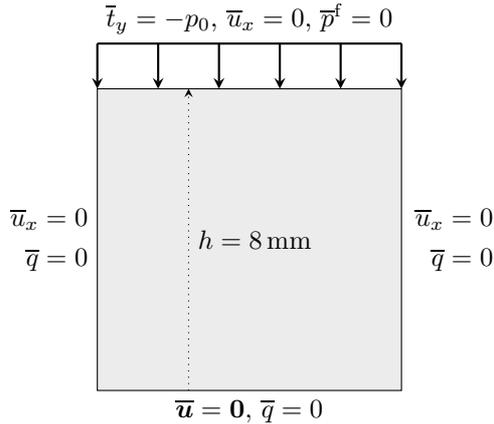

Terzaghi's problem is a classical one-dimensional consolidation problem with an analytical solution, which makes it suitable for code validation. A saturated porous medium subjected an external loading under plane-strain condition is considered where the fluid is allowed to dissipate only at the top boundary, hence resulting in a one-dimensional consolidation. A no flux boundary condition is assumed for the lateral and bottom boundaries. The displacement boundary conditions are such that the lateral sides are constrained from horizontal deformation and the bottom boundary is fixed in both the horizontal and vertical directions. The external load is applied as a Neumann traction $ p_0 $ at the top boundary. The domain and boundary conditions considered are shown in Figure~\ref{fig:terzhagi}.

The analytical solution for the pressure field as a function of time and space is given by:
\begin{align*}
\frac{p^\rmf(t,y)}{p_0}
&=
\frac{4}{\pi} \sum_{i=1}^\infty
\frac{{(-1)}^{i-1}}{2i-1}
\exp \left[ -{(2i-1)}^2 \frac{\pi^2 t_\rms}{4} \right]
\cos \left[ (2i-1) \frac{\pi y}{2h} \right]
\end{align*}
where the dimensionless time $ t_\rms $ is given as a function of the consolidation coefficient $ c_\rmv $ and drainage path $ h $ (total height for one-way drainage) by:

\begin{equation}
t_\rms = \frac{c_\rmv}{h^2}t.
\end{equation}

The consolidation coefficient $ c_\rmv $ is given by:

\begin{equation}
c_\rmv = \frac{(1-\nu)E\kappa}{(1+\nu)(1-2\nu)}
\end{equation}

\begin{table}[t]
    \caption{Terzaghi's problem: Load and material parameters.}
    \label{t:terzaghi_params}
    \begin{center}
        \begin{tabular}{l c r}
            \toprule
            \textbf{Parameter} & \textbf{Value} & \textbf{Unit} \\
            \midrule
            External load, $ p_0 $ & $ 1.0\times10^6 $ & $ \mathrm{Pa} $ \\
            Hydraulic conductivity, $ k $ & $ 1.962\times10^{-14} $ & $ \mathrm{m^2} $ \\
            Biot's coefficient, $ \alpha $ & 1.0 & $ - $ \\
            Young's modulus, $ E $ & $ 6.0\times10^6 $ & $ \mathrm{Pa} $ \\
            Poisson's ratio, $ \nu $ & 0.4 & $ - $ \\
            Storativity, $ c $ & 0 & $ \mathrm{Pa}^{-1} $ \\
            Body forces, $ \bm b $ & $ \bf 0 $ & $ \mathrm{N} $ \\
            \bottomrule 
        \end{tabular}
    \end{center}
\end{table}

The material parameters used for this problem are given in Table~\ref{t:terzaghi_params}, as used in \cite{irzal2013isogeometric}. The choice of the storativity value $c=0$ effectively corresponds to assuming incompressible solid grains and an incompressible fluid.

The Terzaghi verification problem is simulated in a mixed and equal order formulation for comparison. The polynomial degrees considered for the pressure are $ p_\rmp = 1,2,3 $. The corresponding values for the displacement in a mixed formulation are $ p_\rmu = 2,3,4 $. The number of elements used in the simulation is  $N_\rme=\num{72}$. Critical and sub-critical time step sizes are considered to study the sensitivity of the simulations to temporal discretization and to evaluate accuracy of the solution for small time step sizes. The critical time step is calculated according to the relation derived in \cite{vermeer1981accuracy}.

\FloatBarrier

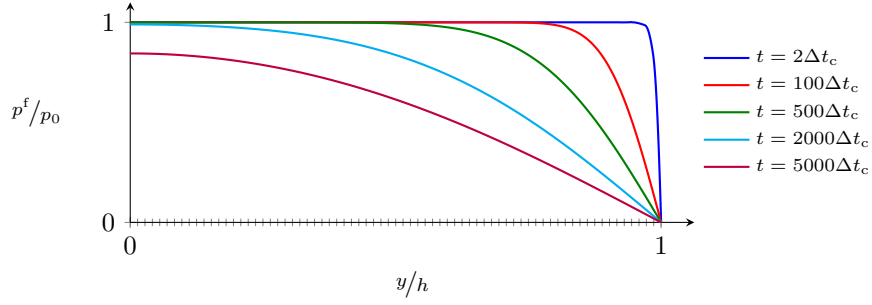
\begin{figure}[t]
    \centering
    \begin{tikzpicture}
    \begin{axis}[
    xmin=0, xmax=0.0085,
    ymin=0, ymax=1.1,
    width=0.6\textwidth,
    height=0.3\textwidth,
    xtick={0,0.008},
    xticklabels={$0$,$1$},
    minor xtick=data,
    ytick={0,1},
    scaled ticks=false,
    axis lines=left,
    xlabel={$\nicefrac{y}{h}$},
    ylabel={$\nicefrac{p^\rmf}{p_0}$},
    xlabel style={at={(axis description cs:0.5,0)}},
    ylabel style={rotate=-90, at={(axis description cs:0,0.5)}},
    legend style={
        at={(axis description cs:1,0.5)},
        anchor=west,
        font=\scriptsize,
        draw=none,
    },
    legend cell align=left,
    ]
    \addplot[draw=none]
    table[x index={0}, y index={1}]{data/terzhagi-crit/ticks.csv};
    \addplot[mark=none, thick, blue]
    table[x index={0}, y index={1}]{data/terzhagi-crit/data-0.csv};
    \addplot[mark=none, thick, red]
    table[x index={0}, y index={1}]{data/terzhagi-crit/data-1.csv};
    \addplot[mark=none, thick, green!50!black]
    table[x index={0}, y index={1}]{data/terzhagi-crit/data-2.csv};
    \addplot[mark=none, thick, cyan]
    table[x index={0}, y index={1}]{data/terzhagi-crit/data-3.csv};
    \addplot[mark=none, thick, purple]
    table[x index={0}, y index={1}]{data/terzhagi-crit/data-4.csv};
    \legend{
        ,
        $t=2\Delta t_\rmc$,
        $t=100\Delta t_\rmc$,
        $t=500\Delta t_\rmc$,
        $t=2000\Delta t_\rmc$,
        $t=5000\Delta t_\rmc$,
    }
    \end{axis}
    \end{tikzpicture}
    \caption{Numerical solution to the Terzhagi problem with $p_\rmp=1$, $p_\rmu=2$ and  and $N_e=72$ using critical time step.    }\label{fig:terzhagi:sol}
\end{figure}   

\begin{figure}[t]
    \begin{subfigure}[b]{1.0\textwidth}
        \centering
        \begin{tikzpicture}
        \begin{axis}[
        xmin=0, xmax=1,
        ymin=0, ymax=1,
        width=1.0\textwidth,
        scale only axis,
        height=2ex,
        axis lines=none,
        legend columns=-1,
        legend style={
            draw=none,
            /tikz/every even column/.append style={column sep=3ex},
            at={(axis description cs:0.5,0.5)},
            anchor=north,
        },
        ]
        \addplot[mark=none, thick, blue] plot coordinates {(0,0) (0,0)};
        \addplot[mark=none, thick, red] plot coordinates {(0,0) (0,0)};
        \addplot[mark=none, thick, cyan] plot coordinates {(2,2) (2,2)};
        \addplot[mark=none, thick, green!50!black] plot coordinates {(2,2) (2,2)};
        \legend{$p_\rmp=1$, $p_\rmp=2$, $p_\rmp=3$}
        \end{axis}
        \end{tikzpicture}
    \end{subfigure}
    \\  
    \begin{subfigure}[b]{0.5\textwidth}
        \centering
        \begin{tikzpicture}
        \begin{axis}[
        xmin=0, xmax=0.0085,
        ymin=0, ymax=1.5,
        width=1.0\textwidth,
        height=0.5\textwidth,
        xtick={0,0.008},
        xticklabels={$0$,$1$},
        minor xtick=data,
        ytick={0,1},
        scaled ticks=false,
        axis lines=left,
        xlabel={$\nicefrac{y}{h}$},
        ylabel={$\nicefrac{p^\rmf}{p_0}$},
        xlabel style={at={(axis description cs:0.5,0)}},
        ylabel style={rotate=-90, at={(axis description cs:0,0.5)}},
        legend style={
            at={(axis description cs:1,0.5)},
            anchor=west,
            font=\scriptsize,
            draw=none,
        },
        legend cell align=left,
        ]
        \addplot[draw=none]
        table[x index={0}, y index={1}]{data/terzhagi-subcrit-mixed/ticks.csv};
        \addplot[mark=none, thick, blue]
        table[x index={0}, y index={1}]{data/terzhagi-subcrit-mixed/data-1.csv};
        \addplot[mark=none, thick, red]
        table[x index={0}, y index={1}]{data/terzhagi-subcrit-mixed/data-2.csv};
        \addplot[mark=none, thick, cyan]
        table[x index={0}, y index={1}]{data/terzhagi-subcrit-mixed/data-3.csv};               
        \end{axis}
        \end{tikzpicture}    
        \caption{Mixed}\label{fig:terzhagi-subcrit:sol}
    \end{subfigure}  
    \begin{subfigure}[b]{0.5\textwidth}
        \centering
        \begin{tikzpicture}
        \begin{axis}[
        xmin=0, xmax=0.0085,
        ymin=0, ymax=1.5,
        width=1.0\textwidth,
        height=0.5\textwidth,
        xtick={0,0.008},
        xticklabels={$0$,$1$},
        minor xtick=data,
        ytick={0,1},
        scaled ticks=false,
        axis lines=left,
        xlabel={$\nicefrac{y}{h}$},
        ylabel={$\nicefrac{p^\rmf}{p_0}$},
        xlabel style={at={(axis description cs:0.5,0)}},
        ylabel style={rotate=-90, at={(axis description cs:0,0.5)}},
        legend style={
            at={(axis description cs:1,0.5)},
            anchor=west,
            font=\scriptsize,
            draw=none,
        },
        legend cell align=left,
        ]
        \addplot[draw=none]
        table[x index={0}, y index={1}]{data/terzhagi-subcrit-equal/ticks.csv};
        \addplot[mark=none, thick, blue]
        table[x index={0}, y index={1}]{data/terzhagi-subcrit-equal/data-1.csv};
        \addplot[mark=none, thick, red]
        table[x index={0}, y index={1}]{data/terzhagi-subcrit-equal/data-2.csv};
        \addplot[mark=none, thick, cyan]
        table[x index={0}, y index={1}]{data/terzhagi-subcrit-equal/data-3.csv};             
        \end{axis}
        \end{tikzpicture}    
        \caption{Equal order}\label{fig:terzhagi-subcrit:sol}
    \end{subfigure}   
    \caption{
        Numerical solution to the Terzaghi problem with $N_\rme=72$ using
        a sub-critical time step of $ ~\Delta t = 0.1\Delta t_\rmc $ for different polynomial degrees. All plots are shown for the first time step.}\label{fig:terzhagi-subcrit:sol}
\end{figure}
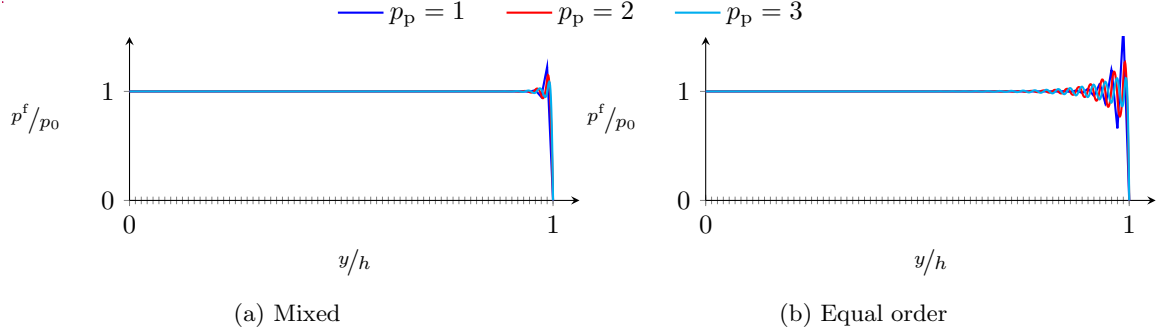

The results from a simulation using the critical time step are shown in Figure~\ref{fig:terzhagi:sol}. A linear solution space is used for the pressure and a quadratic space for the displacement. The results from simulations with a sub-critical time step are shown in Figure~\ref{fig:terzhagi-subcrit:sol} for mixed and equal order cases. The results with the time step size equal to the critical time step show no oscillations in the pressure values. On the other hand, slight oscillations are visible for the sub-critical time step case. These oscillations at very small time steps appear worse for the equal order simulations compared to the mixed simulation. In both cases, the results are observed to improve with increasing polynomial degrees.

\FloatBarrier

\subsection{Terzaghi's Problem: Convergence Study}

Next, a simplified version of the Terzaghi problem is used as a convergence study. We consider a domain with dimensions of $ w \times h = 1 \times 1 $ with the same boundary conditions as in the previous case. For simplicity we choose the following material parameters: $ \alpha=1 $, $ c=0 $, $ E=2/3 $, $ \nu=0.25 $ and $ \kappa=1 $. The external load applied is $ p_0=1 $ and we assume no body forces i.e. $ \bm b = \bm 0 $.

This case was run with an increasing number of degrees of freedom using
polynomial degrees $p_\rmp=1,2,3$ for the pressure and correspondingly $p_\rmu=2,3,4$ for the displacement. In all cases, the time step was kept sufficiently small for the spatial discretization error to dominate and we look at the results at the end of the first time step.

\begin{figure}[h]
    \centering
    \begin{tikzpicture}
    \begin{axis}[
    xmode=log,
    ymode=log,
    width=0.7\textwidth,
    height=0.5\textwidth,
    xlabel={Number of DOFs},
    ylabel={Relative Error, $ \rho_h $},
    legend style={
        at={(axis description cs:0.75,0.75)},
        anchor=west,
        font=\scriptsize,
        draw=none,
    },
    legend cell align=left,
    ]
    \addplot[mark=*, thick, magenta]
    table[x index={0}, y index={1}]{data/conv/lin.csv};
    \addplot[mark=none, dashed, magenta]
    table[x index={0}, y index={1}]{data/conv/comp-lin.csv};
    \addplot[mark=*, thick, blue]
    table[x index={0}, y index={1}]{data/conv/quad.csv};
    \addplot[mark=none, dashed, blue]
    table[x index={0}, y index={1}]{data/conv/comp-quad.csv};
    \addplot[mark=*, thick, red]
    table[x index={0}, y index={1}]{data/conv/cub.csv};
    \addplot[mark=none, dashed, red]
    table[x index={0}, y index={1}]{data/conv/comp-cub.csv};
    \legend{
        $p_\rmp=1$,
        $\mathcal{O}(N^{-1})$,
        $p_\rmp=2$,
        $\mathcal{O}(N^{-\nicefrac{3}{2}})$,
        $p_\rmp=3$,
        $\mathcal{O}(N^{-2})$,
    }
    \end{axis}
    \end{tikzpicture}
    \caption{
        Convergence rates in the relative $L^2$ norm of pressure, for three
        different polynomial degrees.
    } \label{fig:conv}
\end{figure}
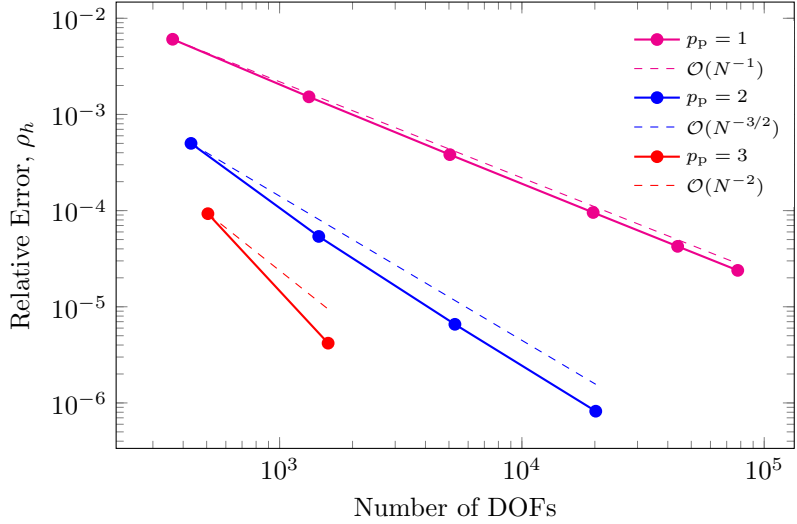

The convergence study is performed by calculating the relative $L^2$ error of the pressure field. The relative error based on the computed pressure values, $ \rho_h $, is calculated from
\begin{equation}
\rho_h = \frac{\|p^\rmf_h - p^\rmf\|_{L^2}}{\|p^\rmf\|_{L^2}}
\end{equation}
where $ p^\rmf_h $ and $ p^\rmf $ are the computed and analytical solution pressures, respectively. The results from the mesh convergence study are shown in Figure~\ref{fig:conv} in terms of plots of the relative error versus the total number of degrees of freedom.  The expected convergence rate based on the analytical solution is also shown. We observe from the results that optimal convergence rates are obtained for all polynomial degrees considered.

\FloatBarrier

\subsection{Low Permeability Layer}

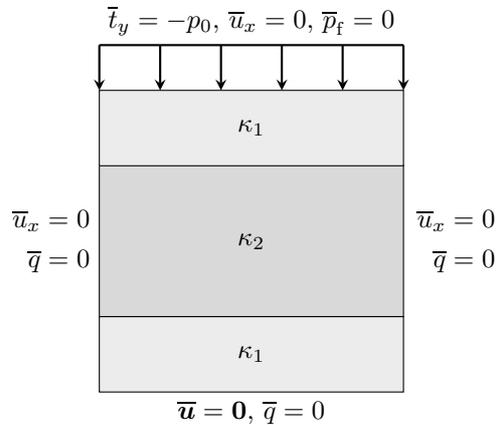
\begin{figure}
    \centering
    \begin{tikzpicture}[scale=4]
    \filldraw[fill=gray!15, draw=black] (0,0) rectangle (1,0.25);
    \filldraw[fill=gray!30, draw=black] (0,0.25) rectangle (1,0.75);
    \filldraw[fill=gray!15, draw=black] (0,0.75) rectangle (1,1);
    \draw [thick] (0,1.15) -- (1,1.15);
    \foreach \x in {0, 0.2, 0.4, 0.6, 0.8, 1}
    \draw[thick, ->,-stealth] (\x,1.15) -- (\x,1);
    \node at (0.5, 0.5) {$\kappa_2$};
    \node at (0.5, 0.125) {$\kappa_1$};
    \node at (0.5, 0.875) {$\kappa_1$};
    \node[anchor=south] at (0.5, 1.15) {
        $\overline{t}_y = -p_0$,
        $\overline{u}_x = 0$,
        $\overline{p}_\rmf = 0$
    };
    \node[anchor=north] at (0.5, 0.0) {
        $\overline{\bm u} = \bm 0$,
        $\overline{q} = 0$
    };
    \node[anchor=west] at (1.0, 0.5) {
        $\begin{aligned}
        \overline{u}_x &= 0 \\
        \overline{q} &= 0
        \end{aligned}$
    };
    \node[anchor=east] at (0.0, 0.5) {
        $\begin{aligned}
        \overline{u}_x &= 0 \\
        \overline{q} &= 0
        \end{aligned}$
    };
    \end{tikzpicture}
    \caption{The Haga problem: Domain and boundary conditions.}\label{fig:haga}
\end{figure}

The next example we consider is the consolidation of a very low permeability layer sandwiched between two high permeability layers, as presented in~\cite{haga2012causes}. A one-dimensional consolidation is assumed by applying the appropriate boundary conditions. The fluid is allowed to dissipate at the top boundary and a no flux condition is defined at the lateral and bottom boundaries. The bottom boundary is fixed from vertical and horizontal displacement and the domain is allowed to deform only in the vertical direction. An external load $p_0$ is applied at the top boundary. The problem setup with the boundary conditions is shown in Figure~\ref{fig:haga}.

The material parameters for this problem are given in Table~\ref{t:haga_params}. Simplified material properties are assumed to focus on the permeability differences of the middle and the bounding layers.

\begin{table}[t]
    \caption{The Haga problem: Load and material parameters.}
    \label{t:haga_params}
    \begin{center}
        \begin{tabular}{l c r}
            \toprule
            \textbf{Parameter} & \textbf{Value} & \textbf{Unit} \\
            \midrule
            External load, $ p_0 $ & $ 1.0 $ & $ \mathrm{Pa} $ \\
            Darcy coefficient, $ k_1/\gamma_\rmf $ & $ 1.0 $ & $ \mathrm{m^2/Pa~s} $ \\
            Darcy coefficient, $ k_2/\gamma_\rmf $ & $ 1.0\times10^{-8} $ & $ \mathrm{m^2/Pa~s} $ \\
            Biot's coefficient, $ \alpha $ & 1.0 & $ - $ \\
            Young's modulus, $ E $ & $ 0.67 $ & $ \mathrm{Pa} $ \\
            Poisson's ratio, $ \nu $ & 0.25 & $ - $ \\
            Storativity, $ c $ & 0 & $ \mathrm{Pa}^{-1} $ \\
            Body forces, $ \bm b $ & $ \bf 0 $ & $ \mathrm{N} $ \\
            \bottomrule 
        \end{tabular}
    \end{center}
\end{table}

The low permeability layer problem is studied using mixed and equal order simulations. The polynomial degrees for the pressure are increased continuously from linear to quartic i.e. $ p_\rmp=1,2,3,4 $. The corresponding polynomial degrees for the displacement in a mixed formulation are $ p_\rmu=2,3,4,5 $. The continuities at the boundaries between the layers are also varied. We consider $C^{0}$ and $C^{p_\rmp-1}$ continuities at these interfaces. In addition, simulations are performed for uniform and graded meshes. The results are presented for these different combinations.

\begin{figure}[t]
    \begin{subfigure}[b]{0.5\textwidth}
        \centering
        \begin{tikzpicture}
        \begin{axis}[
        xmin=0.15, xmax=0.85,
        ymin=0, ymax=2,
        width=1.0\textwidth,
        height=0.5\textwidth,
        xtick={0.25,0.5,0.75},
        minor xtick=data,
        ytick={0,1,2},
        scaled ticks=false,
        axis lines=left,
        xlabel={$\nicefrac{y}{h}$},
        ylabel={$\nicefrac{p^\rmf}{p_0}$},
        xlabel style={at={(axis description cs:0.75,0)}},
        ylabel style={rotate=-90, at={(axis description cs:0,0.75)}},
        ]
        \addplot[draw=none]
        table[x index={0}, y index={1}]{data/haga-unif-max-mixed/ticks.csv};
        \addplot[mark=none, thick, blue]
        table[x index={0}, y index={1}]{data/haga-unif-max-mixed/data-1.csv};
        \addplot[mark=none, thick, red]
        table[x index={0}, y index={1}]{data/haga-unif-max-mixed/data-2.csv};
        \addplot[mark=none, thick, green!50!black]
        table[x index={0}, y index={1}]{data/haga-unif-max-mixed/data-3.csv};
        \addplot[mark=none, thick, cyan]
        table[x index={0}, y index={1}]{data/haga-unif-max-mixed/data-4.csv};
        \end{axis}
        \end{tikzpicture}
    \end{subfigure}
    \begin{subfigure}[b]{0.5\textwidth}
        \centering
        \begin{tikzpicture}
        \begin{axis}[
        xmin=0.15, xmax=0.85,
        ymin=0, ymax=2.0,
        width=1.0\textwidth,
        height=0.5\textwidth,
        xtick={0.25,0.5,0.75},
        minor xtick=data,
        ytick={0,1,2},
        scaled ticks=false,
        axis lines=left,
        xlabel={$\nicefrac{y}{h}$},
        ylabel={$\nicefrac{p^\rmf}{p_0}$},
        xlabel style={at={(axis description cs:0.75,0)}},
        ylabel style={rotate=-90, at={(axis description cs:0,0.75)}},
        ]
        \addplot[draw=none]
        table[x index={0}, y index={1}]{data/haga-unif-max-equal/ticks.csv};
        \addplot[mark=none, thick, blue]
        table[x index={0}, y index={1}]{data/haga-unif-max-equal/data-1.csv};
        \addplot[mark=none, thick, red]
        table[x index={0}, y index={1}]{data/haga-unif-max-equal/data-2.csv};
        \addplot[mark=none, thick, green!50!black]
        table[x index={0}, y index={1}]{data/haga-unif-max-equal/data-3.csv};
        \addplot[mark=none, thick, cyan]
        table[x index={0}, y index={1}]{data/haga-unif-max-equal/data-4.csv};
        \end{axis}
        \end{tikzpicture}
    \end{subfigure}
    \\
    \begin{subfigure}[b]{0.5\textwidth}
        \centering
        \begin{tikzpicture}
        \begin{axis}[
        xmin=0.15, xmax=0.85,
        ymin=0, ymax=2,
        width=1.0\textwidth,
        height=0.5\textwidth,
        xtick={0.25,0.5,0.75},
        minor xtick=data,
        ytick={0,1,2},
        scaled ticks=false,
        axis lines=left,
        xlabel={$\nicefrac{y}{h}$},
        ylabel={$\nicefrac{p^\rmf}{p_0}$},
        xlabel style={at={(axis description cs:0.75,0)}},
        ylabel style={rotate=-90, at={(axis description cs:0,0.75)}},
        ]
        \addplot[draw=none]
        table[x index={0}, y index={1}]{data/haga-unif-min-mixed/ticks.csv};
        \addplot[mark=none, thick, blue]
        table[x index={0}, y index={1}]{data/haga-unif-min-mixed/data-1.csv};
        \addplot[mark=none, thick, red]
        table[x index={0}, y index={1}]{data/haga-unif-min-mixed/data-2.csv};
        \addplot[mark=none, thick, green!50!black]
        table[x index={0}, y index={1}]{data/haga-unif-min-mixed/data-3.csv};
        \addplot[mark=none, thick, cyan]
        table[x index={0}, y index={1}]{data/haga-unif-min-mixed/data-4.csv};
        \end{axis}
        \end{tikzpicture}
    \end{subfigure}
    \begin{subfigure}[b]{0.5\textwidth}
        \centering
        \begin{tikzpicture}
        \begin{axis}[
        xmin=0.15, xmax=0.85,
        ymin=0, ymax=2.0,
        width=1.0\textwidth,
        height=0.5\textwidth,
        xtick={0.25,0.5,0.75},
        minor xtick=data,
        ytick={0,1,2},
        scaled ticks=false,
        axis lines=left,
        xlabel={$\nicefrac{y}{h}$},
        ylabel={$\nicefrac{p^\rmf}{p_0}$},
        xlabel style={at={(axis description cs:0.75,0)}},
        ylabel style={rotate=-90, at={(axis description cs:0,0.75)}},
        ]
        \addplot[draw=none]
        table[x index={0}, y index={1}]{data/haga-unif-min-equal/ticks.csv};
        \addplot[mark=none, thick, blue]
        table[x index={0}, y index={1}]{data/haga-unif-min-equal/data-1.csv};
        \addplot[mark=none, thick, red]
        table[x index={0}, y index={1}]{data/haga-unif-min-equal/data-2.csv};
        \addplot[mark=none, thick, green!50!black]
        table[x index={0}, y index={1}]{data/haga-unif-min-equal/data-3.csv};
        \addplot[mark=none, thick, cyan]
        table[x index={0}, y index={1}]{data/haga-unif-min-equal/data-4.csv};
        \end{axis}
        \end{tikzpicture}
    \end{subfigure}
    \\
    \begin{subfigure}[b]{1.0\textwidth}
        \centering
        \begin{tikzpicture}
        \begin{axis}[
        xmin=0, xmax=1,
        ymin=0, ymax=1,
        width=1.0\textwidth,
        scale only axis,
        height=2ex,
        axis lines=none,
        legend columns=-1,
        legend style={
            draw=none,
            /tikz/every even column/.append style={column sep=3ex},
            at={(axis description cs:0.5,0.5)},
            anchor=north,
        },
        ]
        \addplot[mark=none, thick, blue] plot coordinates {(0,0) (0,0)};
        \addplot[mark=none, thick, red] plot coordinates {(0,0) (0,0)};
        \addplot[mark=none, thick, green!50!black] plot coordinates {(2,2) (2,2)};
        \addplot[mark=none, thick, cyan] plot coordinates {(2,2) (2,2)};
        \legend{$p_\rmp=1$, $p_\rmp=2$, $p_\rmp=3$, $p_\rmp=4$}
        \end{axis}
        \end{tikzpicture}
    \end{subfigure}
    \caption{
        Numerical solution to the Haga problem using $N_\rme=60$ uniform elements and
        $\Delta t=\SI{1}{\second}$. All figures are shown after two time steps. On
        the left the mixed order method, and on the right the equal order method.
        The continuity in the boundary layer is $C^{p_\rmp-1}$ in the top row, and $C^0$
        in the bottom row.
    }\label{fig:haga:unifsol}
\end{figure}

The results from simulations with a uniformly refined mesh are shown in Figure~\ref{fig:haga:unifsol} for the mixed and equal order cases. Severe pressure oscillations are observed within the low permeability layer for the equal order simulations. Due to its high permeability, the fluid in the top layer dissipates very quickly for the time step size considered here i.e. $ \Delta t=1s $. The pressure oscillations start as soon as the fluid in the low permeability layer starts dissipating. The results improve with increasing polynomial degrees but some oscillations are still seen for a quartic solution space for the pressure, $ p_\rmp=4 $. The results with $C^{0}$ continuities at the material interfaces improve slightly better than with $C^{p_\rmp-1}$ continuity since a $C^{0}$ continuity is a more accurate representation of material interfaces. The pressure oscillations in the mixed simulations are less severe and are localized at the boundary between the low permeability and bottom layers. These again decrease with increasing polynomial degrees and a $C^{0}$ continuity at the material interfaces.

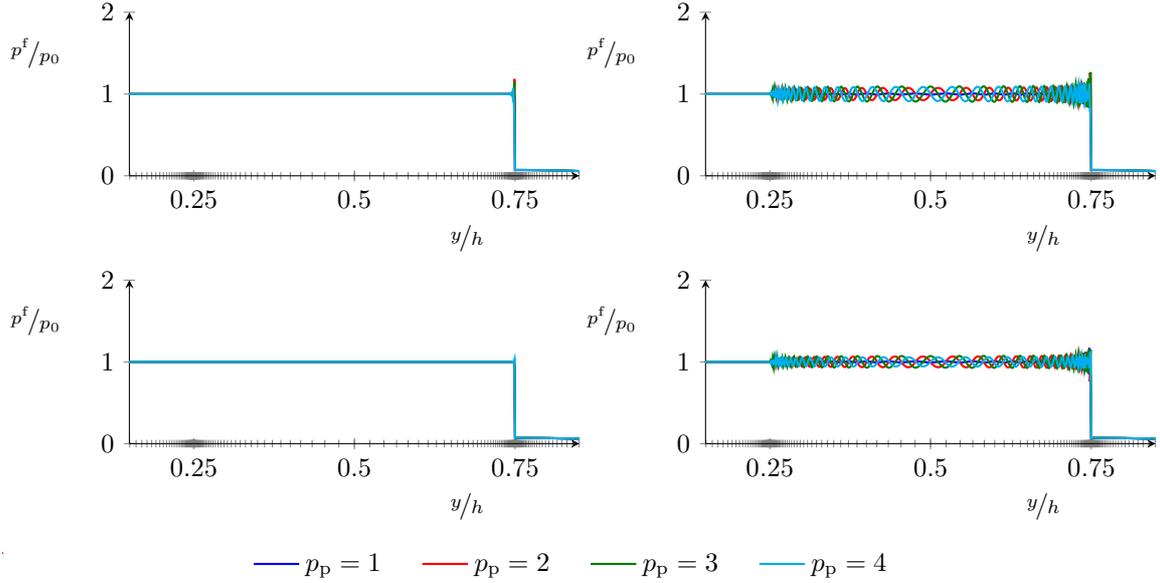
\begin{figure}[h]
    \begin{subfigure}[b]{0.5\textwidth}
        \centering
        \begin{tikzpicture}
        \begin{axis}[
        xmin=0.15, xmax=0.85,
        ymin=0, ymax=2,
        width=1.0\textwidth,
        height=0.5\textwidth,
        xtick={0.25,0.5,0.75},
        minor xtick=data,
        ytick={0,1,2},
        scaled ticks=false,
        axis lines=left,
        xlabel={$\nicefrac{y}{h}$},
        ylabel={$\nicefrac{p^\rmf}{p_0}$},
        xlabel style={at={(axis description cs:0.75,0)}},
        ylabel style={rotate=-90, at={(axis description cs:0,0.75)}},
        ]
        \addplot[draw=none]
        table[x index={0}, y index={1}]{data/haga-grad-max-mixed/ticks.csv};
        \addplot[mark=none, thick, blue]
        table[x index={0}, y index={1}]{data/haga-grad-max-mixed/data-1.csv};
        \addplot[mark=none, thick, red]
        table[x index={0}, y index={1}]{data/haga-grad-max-mixed/data-2.csv};
        \addplot[mark=none, thick, green!50!black]
        table[x index={0}, y index={1}]{data/haga-grad-max-mixed/data-3.csv};
        \addplot[mark=none, thick, cyan]
        table[x index={0}, y index={1}]{data/haga-grad-max-mixed/data-4.csv};
        \end{axis}
        \end{tikzpicture}
    \end{subfigure}
    \begin{subfigure}[b]{0.5\textwidth}
        \centering
        \begin{tikzpicture}
        \begin{axis}[
        xmin=0.15, xmax=0.85,
        ymin=0, ymax=2.0,
        width=1.0\textwidth,
        height=0.5\textwidth,
        xtick={0.25,0.5,0.75},
        minor xtick=data,
        ytick={0,1,2},
        scaled ticks=false,
        axis lines=left,
        xlabel={$\nicefrac{y}{h}$},
        ylabel={$\nicefrac{p^\rmf}{p_0}$},
        xlabel style={at={(axis description cs:0.75,0)}},
        ylabel style={rotate=-90, at={(axis description cs:0,0.75)}},
        ]
        \addplot[draw=none]
        table[x index={0}, y index={1}]{data/haga-grad-max-equal/ticks.csv};
        \addplot[mark=none, thick, blue]
        table[x index={0}, y index={1}]{data/haga-grad-max-equal/data-1.csv};
        \addplot[mark=none, thick, red]
        table[x index={0}, y index={1}]{data/haga-grad-max-equal/data-2.csv};
        \addplot[mark=none, thick, green!50!black]
        table[x index={0}, y index={1}]{data/haga-grad-max-equal/data-3.csv};
        \addplot[mark=none, thick, cyan]
        table[x index={0}, y index={1}]{data/haga-grad-max-equal/data-4.csv};
        \end{axis}
        \end{tikzpicture}
    \end{subfigure}
    \\
    \begin{subfigure}[b]{0.5\textwidth}
        \centering
        \begin{tikzpicture}
        \begin{axis}[
        xmin=0.15, xmax=0.85,
        ymin=0, ymax=2,
        width=1.0\textwidth,
        height=0.5\textwidth,
        xtick={0.25,0.5,0.75},
        minor xtick=data,
        ytick={0,1,2},
        scaled ticks=false,
        axis lines=left,
        xlabel={$\nicefrac{y}{h}$},
        ylabel={$\nicefrac{p^\rmf}{p_0}$},
        xlabel style={at={(axis description cs:0.75,0)}},
        ylabel style={rotate=-90, at={(axis description cs:0,0.75)}},
        ]
        \addplot[draw=none]
        table[x index={0}, y index={1}]{data/haga-grad-min-mixed/ticks.csv};
        \addplot[mark=none, thick, blue]
        table[x index={0}, y index={1}]{data/haga-grad-min-mixed/data-1.csv};
        \addplot[mark=none, thick, red]
        table[x index={0}, y index={1}]{data/haga-grad-min-mixed/data-2.csv};
        \addplot[mark=none, thick, green!50!black]
        table[x index={0}, y index={1}]{data/haga-grad-min-mixed/data-3.csv};
        \addplot[mark=none, thick, cyan]
        table[x index={0}, y index={1}]{data/haga-grad-min-mixed/data-4.csv};
        \end{axis}
        \end{tikzpicture}
    \end{subfigure}
    \begin{subfigure}[b]{0.5\textwidth}
        \centering
        \begin{tikzpicture}
        \begin{axis}[
        xmin=0.15, xmax=0.85,
        ymin=0, ymax=2.0,
        width=1.0\textwidth,
        height=0.5\textwidth,
        xtick={0.25,0.5,0.75},
        minor xtick=data,
        ytick={0,1,2},
        scaled ticks=false,
        axis lines=left,
        xlabel={$\nicefrac{y}{h}$},
        ylabel={$\nicefrac{p^\rmf}{p_0}$},
        xlabel style={at={(axis description cs:0.75,0)}},
        ylabel style={rotate=-90, at={(axis description cs:0,0.75)}}, 
        ]
        \addplot[draw=none]
        table[x index={0}, y index={1}]{data/haga-grad-min-equal/ticks.csv};
        \addplot[mark=none, thick, blue]
        table[x index={0}, y index={1}]{data/haga-grad-min-equal/data-1.csv};
        \addplot[mark=none, thick, red]
        table[x index={0}, y index={1}]{data/haga-grad-min-equal/data-2.csv};
        \addplot[mark=none, thick, green!50!black]
        table[x index={0}, y index={1}]{data/haga-grad-min-equal/data-3.csv};
        \addplot[mark=none, thick, cyan]
        table[x index={0}, y index={1}]{data/haga-grad-min-equal/data-4.csv};
        \end{axis}
        \end{tikzpicture}
    \end{subfigure}
    \\
    \begin{subfigure}[b]{1.0\textwidth}
        \centering
        \begin{tikzpicture}
        \begin{axis}[
        xmin=0, xmax=1,
        ymin=0, ymax=1,
        width=1.0\textwidth,
        scale only axis,
        height=2ex,
        axis lines=none,
        legend columns=-1,
        legend style={
            draw=none,
            /tikz/every even column/.append style={column sep=3ex},
            at={(axis description cs:0.5,0.5)},
            anchor=north,
        },
        ]
        \addplot[mark=none, thick, blue] plot coordinates {(0,0) (0,0)};
        \addplot[mark=none, thick, red] plot coordinates {(0,0) (0,0)};
        \addplot[mark=none, thick, green!50!black] plot coordinates {(2,2) (2,2)};
        \addplot[mark=none, thick, cyan] plot coordinates {(2,2) (2,2)};
        \legend{$p_\rmp=1$, $p_\rmp=2$, $p_\rmp=3$, $p_\rmp=4$}
        \end{axis}
        \end{tikzpicture}
    \end{subfigure}
    \caption{
        Numerical solution to the Haga problem using a graded mesh with small
        elements near the boundary layer and $\Delta t=\SI{1}{\second}$. All figures
        are shown after two time steps. On the left the mixed order method, and on
        the right the equal order method. The continuity in the boundary layer is
        $C^{p_\rmp-1}$ in the top row, and $C^0$ in the bottom row.
    }\label{fig:haga:gradsol}
\end{figure}

Simulations with a graded mesh are also performed for the different combination of polynomial degrees and interface continuities. The graded mesh is generated such that more elements are concentrated at the material interfaces. The results from this case are shown in Figure~\ref{fig:haga:gradsol}. The pressure oscillations in the equal order case improve significantly in this case compared to the results from uniform mesh refinement. However, the oscillations still occur throughout the low permeability layer. The equal order results for linear basis functions show a slightly strange behavior in that the oscillations are lesser within the low permeability layer than for higher order elements, but show slightly higher oscillations at the top material interface. The results are again better with a $C^{0}$ continuity at the material interfaces. The mixed simulation results also improve with a graded mesh. Almost no oscillations are noticed for combinations of higher polynomial degrees and $C^{0}$ continuities at the material interfaces.

\section{Conclusions}
\label{sec:conclusions}

Mixed isogeometric analysis of poroelasticity is presented where different order of polynomials are used for the displacement and pore pressure field variables. Numerical studies on Terzaghi's classical one-dimensional consolidation problem and consolidation of a layered soil with a middle low permeability layer are presented. The results from mixed polynomial order simulations are compared with equal order analyses. For Terzaghi's one-dimensional consolidation problem, the pore pressure oscillations are investigated when a time step size less than the critical value is used. The oscillations were observed to be higher in the equal order simulations compared to the mixed order results. The oscillations are not completely removed in the mixed isogeometric simulations but it is observed that the they tend to decrease with increasing polynomial orders for the pore pressure. This is illustrated by the convergence of the relative $ L^2 $ norm of the pore pressure error for varying polynomial orders. The low permeability layer problem showed similar trends in the pore pressure oscillations i.e. the equal order simulations resulted in worse pore pressure oscillations compared to the mixed results. Again, in both cases, the oscillations decreased with increasing polynomial orders. The use of a graded mesh, where the knot spans are concentrated at the interfaces between the low permeability and other layers, resulted in much lower oscillations both in the equal order and mixed cases. This indicates the potential of adaptive refinement for such class of problems.    

\section*{Acknowledgement}

This work is financially supported by the Research Council of Norway and industrial partners through the research project SAMCoT, Sustainable Arctic Marine and Coastal Technology. The authors gratefully acknowledge the support.

\bibliographystyle{plain}
\bibliography{References}

\begin{thebibliography}{10}

\bibitem{bekele2015adaptive}
Yared~W Bekele, Trond Kvamsdal, Arne~M Kvarving, and Steinar Nordal.
\newblock Adaptive isogeometric finite element analysis of steady-state
  groundwater flow.
\newblock {\em International Journal for Numerical and Analytical Methods in
  Geomechanics}, 2015.

\bibitem{biot1955theory}
M~Av Biot.
\newblock Theory of elasticity and consolidation for a porous anisotropic
  solid.
\newblock {\em Journal of Applied Physics}, 26(2):182--185, 1955.

\bibitem{biot1956general}
MA~Biot.
\newblock General solutions of the equations of elasticity and consolidation
  for a porous material.
\newblock {\em J. appl. Mech}, 23(1):91--96, 1956.

\bibitem{biot1941general}
Maurice~A Biot.
\newblock General theory of three-dimensional consolidation.
\newblock {\em Journal of applied physics}, 12(2):155--164, 1941.

\bibitem{booker1975investigation}
John~R Booker and JC~Small.
\newblock An investigation of the stability of numerical solutions of {B}iot's
  equations of consolidation.
\newblock {\em International Journal of Solids and Structures}, 11(7):907--917,
  1975.

\bibitem{cottrell2009iga}
J.~Austin Cottrell, Thomas J.~R. Hughes, and Yuri Bazilevs.
\newblock {\em Isogeometric Analysis : Toward Integration of {CAD} and {FEA}}.
\newblock Wiley, Chichester, West Sussex, U.K., Hoboken, NJ, 2009.

\bibitem{de1988historical}
R~de~Boer and W~Ehlers.
\newblock A historical review of the formulation of porous media theories.
\newblock {\em Acta Mechanica}, 74(1-4):1--8, 1988.

\bibitem{dureisseix2003latin}
David Dureisseix, Pierre Ladev{\`e}ze, and Bernard~A Schrefler.
\newblock A {LATIN} computational strategy for multiphysics problems:
  application to poroelasticity.
\newblock {\em International Journal for Numerical Methods in Engineering},
  56(10):1489--1510, 2003.

\bibitem{hamalawi1999posteriori}
Ashraf El-Hamalawi and MD~Bolton.
\newblock An a posteriori error estimator for plane-strain geotechnical
  analyses.
\newblock {\em Finite elements in analysis and design}, 33(4):335--354, 1999.

\bibitem{hamalawi2002posteriori}
Ashraf El-Hamalawi and MD~Bolton.
\newblock A-posteriori error estimation in axisymmetric geotechnical analyses.
\newblock {\em Computers and Geotechnics}, 29(8):587--607, 2002.

\bibitem{ferronato2001ill}
Massimiliano Ferronato, Giuseppe Gambolati, and Pietro Teatini.
\newblock Ill-conditioning of finite element poroelasticity equations.
\newblock {\em International Journal of Solids and Structures},
  38(34):5995--6014, 2001.

\bibitem{gambolati2001numerical}
Giuseppe Gambolati, Giorgio Pini, and Massimiliano Ferronato.
\newblock Numerical performance of projection methods in finite element
  consolidation models.
\newblock {\em International Journal for Numerical and Analytical Methods in
  Geomechanics}, 25(14):1429--1447, 2001.

\bibitem{ghaboussi1973flow}
Jamshid Ghaboussi and Edward~L Wilson.
\newblock Flow of compressible fluid in porous elastic media.
\newblock {\em International Journal for Numerical Methods in Engineering},
  5(3):419--442, 1973.

\bibitem{haga2012causes}
Joachim~Berdal Haga, Harald Osnes, and Hans~Petter Langtangen.
\newblock On the causes of pressure oscillations in low-permeable and
  low-compressible porous media.
\newblock {\em International Journal for Numerical and Analytical Methods in
  Geomechanics}, 36(12):1507--1522, 2012.

\bibitem{hughes2005iga}
T.~J.~R. Hughes, J.~A. Cottrell, and Y.~Bazilevs.
\newblock Isogeometric analysis: {CAD}, finite elements, {NURBS}, exact
  geometry and mesh refinement.
\newblock {\em Computer Methods in Applied Mechanics and Engineering},
  194(39–41):4135--4195, 2005.

\bibitem{hwang1971solutions}
CoT Hwang, NR~Morgenstern, and DW~Murray.
\newblock On solutions of plane strain consolidation problems by finite element
  methods.
\newblock {\em Canadian Geotechnical Journal}, 8(1):109--118, 1971.

\bibitem{irzal2013isogeometric}
Faisal Irzal, Joris~JC Remmers, Clemens~V Verhoosel, and Ren{\'e} de~Borst.
\newblock Isogeometric finite element analysis of poroelasticity.
\newblock {\em International Journal for Numerical and Analytical Methods in
  Geomechanics}, 37(12):1891--1907, 2013.

\bibitem{korsawe2006finite}
Johannes Korsawe, Gerhard Starke, Wenqing Wang, and Olaf Kolditz.
\newblock Finite element analysis of poro-elastic consolidation in porous
  media: Standard and mixed approaches.
\newblock {\em Computer Methods in Applied Mechanics and Engineering},
  195(9):1096--1115, 2006.

\bibitem{larsson2015sequential}
Fredrik Larsson and Kenneth Runesson.
\newblock A sequential-adaptive strategy in space-time with application to
  consolidation of porous media.
\newblock {\em Computer Methods in Applied Mechanics and Engineering},
  288:146--171, 2015.

\bibitem{lewis1998finite}
R.W. Lewis and B.A. Schrefler.
\newblock {\em The finite element method in the static and dynamic deformation
  and consolidation of porous media}.
\newblock Numerical methods in engineering. John Wiley, 1998.

\bibitem{murad1992improved}
M{\'a}rcio~A Murad and Abimael~FD Loula.
\newblock Improved accuracy in finite element analysis of {B}iot's
  consolidation problem.
\newblock {\em Computer Methods in Applied Mechanics and Engineering},
  95(3):359--382, 1992.

\bibitem{murad1994stability}
M{\'a}rcio~A Murad and Abimael~FD Loula.
\newblock On stability and convergence of finite element approximations of
  {B}iot's consolidation problem.
\newblock {\em International Journal for Numerical Methods in Engineering},
  37(4):645--667, 1994.

\bibitem{phillips2007coupling1}
Phillip~Joseph Phillips and Mary~F Wheeler.
\newblock A coupling of mixed and continuous {G}alerkin finite element methods
  for poroelasticity {I}: the continuous-in-time case.
\newblock {\em Computational Geosciences}, 11(2):131--144, 2007.

\bibitem{phillips2007coupling2}
Phillip~Joseph Phillips and Mary~F Wheeler.
\newblock A coupling of mixed and continuous {G}alerkin finite element methods
  for poroelasticity {II}: the discrete-in-time case.
\newblock {\em Computational Geosciences}, 11(2):145--158, 2007.

\bibitem{phillips2008coupling3}
Phillip~Joseph Phillips and Mary~F Wheeler.
\newblock A coupling of mixed and discontinuous {G}alerkin finite-element
  methods for poroelasticity.
\newblock {\em Computational Geosciences}, 12(4):417--435, 2008.

\bibitem{reed1984investigation}
MB~Reed.
\newblock An investigation of numerical errors in the analysis of consolidation
  by finite elements.
\newblock {\em International journal for numerical and analytical methods in
  geomechanics}, 8(3):243--257, 1984.

\bibitem{sandhu1985special}
Ranbir~S Sandhu, Shyan~Chyun Lee, and Hwie-Ing The.
\newblock Special finite elements for analysis of soil consolidation.
\newblock {\em International journal for numerical and analytical methods in
  geomechanics}, 9(2):125--147, 1985.

\bibitem{sandhu1977numerical}
Ranbir~S Sandhu, Honho Liu, and Kamar~J Singh.
\newblock Numerical performance of some finite element schemes for analysis of
  seepage in porous elastic media.
\newblock {\em International Journal for Numerical and Analytical Methods in
  Geomechanics}, 1(2):177--194, 1977.

\bibitem{sandhu1969finite}
Ranbir~S Sandhu and Edward~L Wilson.
\newblock Finite-element analysis of seepage in elastic media.
\newblock {\em Journal of the Engineering Mechanics Division}, 95(3):641--652,
  1969.

\bibitem{sloan1999biot1}
Scott~W Sloan and Andrew~J Abbo.
\newblock Biot consolidation analysis with automatic time stepping and error
  control part 1: theory and implementation.
\newblock {\em International Journal for Numerical and Analytical Methods in
  Geomechanics}, 23(6):467--492, 1999.

\bibitem{sloan1999biot2}
Scott~W Sloan and Andrew~J Abbo.
\newblock Biot consolidation analysis with automatic time stepping and error
  control: Part 2: Applications.
\newblock {\em International Journal for Numerical and Analytical Methods in
  Geomechanics}, 23(6):493--529, 1999.

\bibitem{tchonkova2008new}
Maria Tchonkova, John Peters, and Stein Sture.
\newblock A new mixed finite element method for poro-elasticity.
\newblock {\em International journal for numerical and analytical methods in
  geomechanics}, 32(6):579--606, 2008.

\bibitem{terzaghi1923berechnung}
K~von Terzaghi.
\newblock Die berechnung der durchlassigkeitsziffer des tones aus dem verlauf
  der hydrodynamischen spannungserscheinungen.
\newblock {\em Sitzungsberichte der Akademie der Wissenschaften in Wien,
  Mathematisch-Naturwissenschaftliche Klasse, Abteilung IIa}, 132:125--138,
  1923.

\bibitem{terzaghi1925erdbaumechanik}
Karl Terzaghi et~al.
\newblock Erdbaumechanik auf bodenphysikalischer grundlage.
\newblock 1925.

\bibitem{vermeer1981accuracy}
PA~Vermeer and A~Verruijt.
\newblock An accuracy condition for consolidation by finite elements.
\newblock {\em International Journal for numerical and analytical methods in
  geomechanics}, 5(1):1--14, 1981.

\bibitem{zhu2004numerical}
Guofu Zhu, Jian-Hua Yin, and Shun-tim Luk.
\newblock Numerical characteristics of a simple finite element formulation for
  consolidation analysis.
\newblock {\em Communications in Numerical Methods in Engineering},
  20(10):767--775, 2004.

\end{thebibliography}

\end{document}